\documentclass{article}
\usepackage[a4paper,left=2cm,right=2cm,top=2cm,bottom=2cm]{geometry}
\usepackage[cp1250]{inputenc}
\usepackage{amsfonts,amsmath,upref,amssymb,amsthm,amsxtra,amscd,enumerate}
\usepackage{tikz}
\DeclareMathOperator{\diver}{div}

\DeclareMathOperator{\supp}{supp}

\newcommand{\vv}{{\bf v}}
\newcommand{\vm}{{\bf m}}

\newcommand{\pat}{\partial_t}
\newcommand{\Div}{\diver_x}

\newtheorem{theorem}{Theorem}[section]

\newtheorem{definition}[theorem]{Definition}

\numberwithin{equation}{section}

\title{A Note on Measure-Valued Solutions to the Full Euler System}
\date{}
\author{V\' aclav M\' acha\footnote{Institute of Mathematics of the Academy of Sciences of the Czech Republic,
e-mail: macha@math.cas.cz} 
\and Emil Wiedemann\footnote{Institute of Applied Analysis, Ulm University, Germany, e-mail: emil.wiedemann@uni-ulm.de}}
\begin{document}
\maketitle

\begin{abstract}
We construct two particular solutions of the full Euler system which emanate from the same initial data. Our aim is to show that the convex combination of these two solutions form a measure-valued solution which may not be approximated by a sequence of weak solutions. As a result, the weak* closure of the set of all weak solutions, considered as parametrized measures, is not equal to the space of all measure-valued solutions. This is in stark contrast with the incompressible Euler equations. 
\end{abstract}

\section{Introduction}
In the context of fluid dynamics, measure-valued solutions were first studied by DiPerna and Majda~\cite{dipernamajda}, who developed an appropriate mathematical framework and showed existence for such solutions of the incompressible Euler equations. Measure-valued solutions describe the one-point statistics of a fluid, i.e., they give the probability distribution of the fluid velocity (and other state variables like density or temperature) at a given point in time and space. If one is willing to accept such a probabilistic description rather than a deterministic one (which would, of course, contain more information), then one easily obtains a solution for any initial data, bypassing the notorious problem of non-interchangeability of weak limits and nonlinearities.

Measure-valued solutions are sometimes thought of as a ``cheap way out'' of the fundamental lack of compactness for inviscid fluid models, and are criticized as not containing enough interesting information. Yet, in recent years, the concept has been intensely studied again, as it turned out to have several merits after all: First, measure-valued solutions, despite representing a very weak notion of solution, enjoy a weak-strong stability property that entails important consequences for singular limits, numerical approximation, and long-time behaviour; this property is known for many systems of fluid mechanics, including the incompressible~\cite{brenierDLSz}, isentropic compressible~\cite{GwSwWi}, and full compressible~\cite{BrFe} Euler systems, the isentropic compressible Navier-Stokes equations~\cite{FeGwSwWi}, and the Navier-Stokes-Fourier equations~\cite{BrFeNo}, see also the survey~\cite{wsusurvey}. One should remark that weak-strong uniqueness holds only for \emph{admissible} measure-valued solutions, which comply with appropriate energy or entropy inequalities. 

Motivated by numerical simulations, Fjordholm et al.~\cite{FjMiTa} argued that measure-valued solutions provide for a more suitable notion of solution than weak (distributional) solutions; indeed, the numerical computation of unstable shear flows with randomly perturbed initial data yields highly unpredictable results on the level of weak solutions, but apparently stable and regular behavior on the measure-valued level. This phenomenon is of course very plausible in the light of phenomenological turbulence theory~\cite{Frisch}.   

Measure-valued solutions seem like a vast generalization of weak solutions, but are they really? Surprisingly, the answer is `no' for the incompressible Euler equations~\cite{SzWi12}. Indeed, any measure-valued solution is weakly approximated by a sequence of weak solutions (and if the measure is admissible, then the approximating sequence can also be chosen to consist of admissible weak solutions), or in other words: The set of Dirac parametrized measures is weakly* dense in the set of all measure-valued solutions. One might thus say that the notion of measure-valued solution is just a topological closure, but not a substantial extension, of the more classical concept of weak solution.

Looking at~\cite{SzWi12} from a different angle, one could view the result as an instance of a characterization of Young measures generated by sequences with specific properties (viz., being a solution of the incompressible Euler equations). The most classical result of this kind is the characterization of gradient Young measures~\cite{KiPe91, KiPe94}, where \emph{not} every Young measure whose barycenter is a gradient is itself generated by a sequence of gradients. This already indicates that the situation for incompressible Euler is rather unusual. 

In fact, the fact that every measure-valued solution of the incompressible Euler equations is generable is related to the \emph{wave cone} for the corresponding linear constraint; indeed, the wave cone in this case is the whole space. It turns out that this is no longer the case for compressible models. In~\cite{ChFeKrWi}, the wave cone for the isentropic compressible Euler system is determined (and it is not the whole space), and a preliminary application to the generability of measure-valued solutions is given. This result was recently extended in~\cite{GaWi} to yield an admissible measure-valued solution with atomic initial data that is not recovered from a sequence of weak solutions. The result relies on a refined rigidity lemma (see Theorem~\ref{thm.2.1} below) and the construction of `wild' solutions for certain Riemann data in~\cite{ChDLKr}. It is shown, in addition, that the constructed measure-valued solution does not arise as a vanishing viscosity limit, and it is argued that such measure-valued solutions, therefore, should be discarded as unphysical.

The aim of the present contribution is to extend the result from~\cite{GaWi} to the full Euler system. Our main result (Theorem~\ref{mainthm} below) states that there exists a measure-valued solution to the full Euler system with non-constant entropy which is not generated by any sequence of weak solutions. 

This contribution follows a similar strategy as~\cite{GaWi}, but it relies on the convex integration construction from~\cite{KlKrMaMa} for the full Euler system rather than the isentropic construction from~\cite{ChDLKr}. We have tried to keep the presentation as concise as possible; this means that we only briefly recall the relevant facts from~\cite{GaWi,KlKrMaMa} with no extensive discussion. Moreover, our result could be extended in various straightforward ways to include the admissibility condition (which is actually satisfied in our construction, since the scheme in~\cite{KlKrMaMa} produces entropy solutions), concentration measures (which we could simply take to be zero, whence any hypothetical generating sequence would automatically have equi-integrable nonlinearities, cf.~\cite[proof of Theorem 4.13]{GaWi}), or viscosity limits. We have chosen, however, to keep this note as short and simple as possible.    

\section{Preliminaries}
Let  $T>0$. 
We consider the following system on time-space $[0,T]\times\mathbb R^2$:
\begin{equation}\label{f.euler}
\begin{split}
\pat \rho  + \Div (\rho \vv)  &= 0\\
\pat (\rho \vv) + \Div (\rho \vv \otimes \vv) + \nabla p(\varrho,\theta) & = 0\\
\pat \left(\frac 12 \rho |\vv|^2 + \rho e(\rho,\theta)\right) + \Div \left(\left(\frac 12 \rho |\vv|^2 + \rho e(\rho,\theta)\right)\vv\right) + \Div (p(\rho,\theta) \vv) & = 0,
\end{split}
\end{equation}
with unknowns $\rho:[0,T]\times\Omega\to \mathbb R^+_0$, $\vv=(v,u):[0,T]\times \Omega \to \mathbb R^2$ and $\theta:[0,T]\times \Omega \to \mathbb R^+_0$. The functions $e$ and $p$ are interrelated through the Gibbs law which gives rise also to an entropy -- a function $s(\rho,\theta)$ such that 
\begin{equation*}
\theta D s(\rho,\theta) = De(\rho,\theta) + p(\rho,\theta) D\left(\frac 1\rho\right),
\end{equation*}
where $D$ stands for the gradient with respect to $\rho$ and $\theta$. Throughout this paper we consider an ideal gas, i.e., 
$$
p(\rho,\theta) = \rho \theta,\ e(\rho,\theta) = c_v\theta,\ s(\rho,\theta) = \log \left(\frac{\theta^{c_v}}\rho\right),\ \mbox{for some }c_v>0
$$
and, for simplicity, we assume $c_v = 1$. We rewrite \eqref{f.euler} into conservative variables $\rho, \vm=\rho\vv,\ E=\frac 12 \rho |\vv|^2 + \rho e(\rho,\theta)$. Our choice of state variables gives $p = E - \frac 12 \frac{\vm^2}\rho$ and thus we get
\begin{equation}\label{f.e.cons}
\begin{split}
\pat \rho + \Div \vm & = 0\\
\pat \vm + \Div \left(\frac{\vm\otimes \vm}\rho + \mathbb I \left(E -\frac 12 \frac{\vm^2}\rho\right)\right) & = 0\\
\pat E + \diver \left(\left(2E-\frac 12 \frac{\vm^2}\rho\right)\frac \vm\rho\right) & = 0,
\end{split}
\end{equation}
where $\mathbb I$ denotes the $(2\times 2)$ identity matrix.
Furthermore, we may rewrite it as a linear differential system
\begin{equation}
\begin{split}\label{f.e.linear}
\pat \rho + \Div \vm & = 0\\
\pat \vm + \Div U + \nabla E & = 0\\
\pat E + \Div r & = 0
\end{split}
\end{equation}
with the following constraints:
\begin{equation}
\begin{split}\label{constrain}
U &= \frac{\vm\otimes \vm}\rho - \frac 12 \frac{\vm^2}{\rho} \mathbb I\\
r & = \left(\left(2E-\frac 12 \frac{\vm^2}\rho\right)\frac \vm\rho\right).
\end{split}
\end{equation}

We recall that, according to the definition, $U$ is a traceless symetric matrix and thus the system \eqref{f.e.linear} may be rewritten as
\begin{equation}
\diver \left(
\begin{matrix}
\rho & m_1 & m_2\\
m_1 & U_{11} + E & U_{12}\\
m_2 & U_{12} & -U_{11} + E\\
E & r_1 & r_2
\end{matrix}
\right) = 0.
\end{equation}
Here $\diver$ stands for a divergence over the time-space variables $(t,x)$.

\begin{definition}\label{def:mvs}
We say that a family of probability measures $\nu:=\nu_{t,x}\in L^\infty_{w*}([0,T]\times\mathbb R^2, \mathcal P(\mathbb R^+_0\times \mathbb R^2\times \mathbb R^+_0))$ is a measure-valued solution to \eqref{f.euler} with initial data $(\rho^0,\vm^0,E^0)$ if
\begin{itemize}
\item \begin{equation}\int_0^T\int_{\mathbb R^2} \rho(t,x) \pat \varphi(t,x)  + \vm(t,x)\cdot\nabla \varphi(t,x)\ {\rm d}x{\rm d}t +\int_{\mathbb R^2} \rho^0(x) \varphi(0,x){\rm d}x= 0
\end{equation}
for all $\varphi \in C^\infty_c  ([0,T)\times \mathbb R^2)$,
\item \begin{equation}\int_0^T\int_{\mathbb R^2} \vm(t,x)\cdot\pat \varphi(t,x) + \left\langle\nu_{t,x}, \frac{\xi_m\otimes \xi_m}{\xi_\rho} - \mathbb I\left(\xi_E - \frac 12 \frac{\xi_m^2}{\xi_\rho}\right) \right\rangle :\nabla \varphi(t,x)\ {\rm d}x{\rm d}t +\int_{\mathbb R^2} \vm^0(x) \cdot\varphi(0,x){\rm d}x= 0
\end{equation}
for all $\varphi \in C^\infty_c ([0,T)\times \mathbb R^2)^2$,
\item \begin{equation}\int_0^T\int_{\mathbb R^2} E(t,x)\pat \varphi (t,x) + \left\langle\nu_{t,x}, \left(2\xi_E - \frac 12 \frac{\xi_m^2}{\xi_\rho}\right)\frac{\xi_m}{\xi_\rho}\right\rangle \cdot\nabla \varphi(t,x){\rm d}x{\rm d}t+\int_{\mathbb R^2} E^0(x) \varphi (0,x){\rm d}x=0 \end{equation}
for all $\varphi \in C^\infty_c([0,T)\times \mathbb R^2)$.
\end{itemize} 
Here $\xi_\rho\in \mathbb R^+_0$, $\xi_m\in \mathbb R^2$, and $\xi_E\in \mathbb R^+_0$ are dummy variables for $\rho$, $m$ and $E$, meaning that 
$$
\rho(t,x) = \int_0^\infty \xi_\rho \ {\rm d}\nu_{t,x}, \ \vm(t,x) = \int_{\mathbb R^2} \xi_m \ {\rm d}\nu_{t,x},\ E = \int_0^\infty \xi_E \ {\rm d}\nu_{t,x}
$$
and we use the notation
$$
\langle \nu_{t,x},f(\xi_\rho,\xi_m, \xi_E)\rangle =  \int_{\mathbb R^+_0\times \mathbb R^2\times \mathbb R^+_0} f(\xi_\rho,\xi_m, \xi_E)\ {\rm d}\nu_{t,x}.
$$
\end{definition}

For the existence of a measure-valued solution, one needs an extended notion including concentration effects; if these are taken into account, then the existence of measure-valued solutions is known owing to J.~B\v rezina \cite{brezina}. His definition of measure-valued solution is also slightly different from ours in that he used the renormalized entropy balance and the total energy balance instead of the energy balance. In any case, every weak solution $(\rho,\vm,E)$ is also a measure valued solution (in the sense of Definition~\ref{def:mvs}) with $\nu_{t,x} = \delta_{\rho(t,x)}\otimes \delta_{\vm(t,x)}\otimes \delta_{E(t,x)}$. Thus, the existence of infinitely many weak solutions for certain initial data exhibited in \cite{AlKlKrMaMa} and in \cite{KlKrMaMa} shows a fortiori the existence of non-unique measure-valued solutions for these data.

As an immediate consequence of the definition we obtain that every convex combination of two measure valued solutions is also a measure-valued solution, that is, if $\nu, \mu$ are two measure-valued solutions, then so is $\lambda\nu+(1-\lambda)\mu$ for any $\lambda\in [0,1]$.

Our aim is to prove that there exists a measure valued solution to \eqref{f.euler} which cannot be obtained as a limit of weak solutions. To make this precise, we say that a sequence $(z_n):\Omega\to\mathbb R^d$ of measurable functions generates the parametrized measure $(\nu_x)_{x\in\Omega}$ if
\begin{equation*}
f(z_n)\rightharpoonup \langle \nu, f\rangle\quad \text{weakly in $L^1(\Omega)$}
\end{equation*}
for all continuous functions $f:\mathbb R^d\to\mathbb R$ for which $(f(z_n))$ is equi-integrable.

We will take advantage of the following theorem proved in \cite{ChFeKrWi}:
\begin{theorem}\label{thm:constants}
Let $\Omega \subset \mathbb R^N$ be a Lipschitz and bounded domain, $1\leq p<\infty$, and $\mathcal A$ a linear operator of the form
\begin{equation}
\mathcal A z := \sum_{i=1}^N A^{(i)} \frac{\partial z}{\partial x_i},
\end{equation}
where $A^{(i)}$ are $l\times d$ matrices and $z:\mathbb R^N\mapsto \mathbb R^d$ a vector valued function. Let $p\in (1,\infty)$ and $z_1, z_2\in \mathbb R^d$, $z_1\neq z_2$ be two constant states such that 
$$
z_2 - z_1\notin \Lambda,
$$
where $\Lambda$ denotes the wave cone, defined by
\begin{equation*}
\Lambda=\left\{\bar z\in\mathbb R^d: \text{there exists $\xi\in\mathbb R^N\setminus\{0\}$ such that $\mathcal A(\bar z h(\cdot\cdot\xi))=0$ for all $h:\mathbb R\to\mathbb R$}\right\}.
\end{equation*}
 Let further $z_n:\Omega \to \mathbb R^d$ is an equi-integrable family of functions with 
\begin{equation*}
\begin{split}
\|z_n\|_{L^p}& \leq c\\
\mathcal A z_n \to 0\ \mbox{in }W^{-1,r}(\Omega)
\end{split}
\end{equation*}
for some $r\in \left(1,\frac{N}{N-1}\right)$, and assume that $(z_n)$ generates a compactly supported Young measure such that 
\begin{equation*}
\supp [\nu_x]\subset \{\lambda z_1 + (1-\lambda)z_2,\ \lambda \in [0,1]\}\ \mbox{for a.a. }x\in \Omega.
\end{equation*}
Then there exists $z_\infty\in \mathbb R^d$ such that 
$$
z_n\to z_\infty \ \mbox{in }L^p(\Omega).
$$
\end{theorem}

In our setting (i.e., the state vector is $(\rho,m_1,m_2,U_{11}, U_{12}, E, r_1, r_2)$, we have
$$
A^{(1)} = \left(
\begin{matrix}
1&0&0&0&0&0&0&0\\
0&1&0&0&0&0&0&0\\
0&0&1&0&0&0&0&0\\
0&0&0&0&0&1&0&0
\end{matrix}
\right),
$$
$$
A^{(2)} = \left(
\begin{matrix}
0&1&0&0&0&0&0&0\\
0&0&0&1&0&1&0&0\\
0&0&0&0&1&0&0&0\\
0&0&0&0&0&0&1&0
\end{matrix}
\right),
$$
$$
A^{(3)} = \left(
\begin{matrix}
0&0&1&0&0&0&0&0\\
0&0&0&0&1&0&0&0\\
0&0&0&-1&0&1&0&0\\
0&0&0&0&0&0&0&1
\end{matrix}
\right),
$$

For a given point $z = (\rho,m_1,m_2,U_{11}, U_{12}, E, r_1, r_2)$ define a matrix $(Z_{\mathcal A})_{ji}$ as follows:
$$
(Z_{\mathcal A})_{ji} = \sum_{k=1}^d A^{(i)}_{jk} z_k,\ \qquad j=1,\ldots,4,\ i=1,\ldots,3.
$$

As observed in \cite[Section 3.2]{ChFeKrWi}, $z\in \Lambda$ if and only if the corresponding $Z_{\mathcal A}$ satisfies $\mbox{rank}\ Z_{\mathcal A} < 3$. Thus it is enough to take
$$
z_1 = \left(1,1,0,\frac 12, 0,\frac 32, \frac 52,0\right),\ z_2 = \left(\gamma, 1,0,\frac 1{2\gamma} ,0, \frac 3{2\gamma}, \frac 5{2\gamma^2}, 0 \right).
$$ 
Trivially, both $z_1$ and $z_2$ are solutions to \eqref{f.e.cons} since they are constant. Moreover, there exists $\gamma$ such that $z_1-z_2\notin \Lambda$. Indeed, the corresponding $Z_{\mathcal A}$ is of the form 
$$
\left(\begin{matrix}
1-\gamma & 0 & 0 \\
0 &2\left(1-\frac 1\gamma\right) & 0 \\
0 & 0 &  1 - \frac 1\gamma \\
\frac32(1-\frac 1{\gamma}) & \frac52(1-\frac 1{\gamma^2}) & 0
\end{matrix}\right)
$$
and the determinant of the $3\times 3$ submatrix $\left(Z_{\mathcal A}\right)_{i,j = 1}^3$ is $2(1-\gamma)(1-\frac 1\gamma)^2$ and thus is nonzero for all $\gamma \neq 1$. According to Theorem \ref{thm:constants}, the measure valued solution $\nu_{t,x} = \frac 12\delta_{z_1} + \frac 12 \delta_{z_2}$ may not be approximated by a sequence of weak solutions. This is the cheapest way how to produce a solution of the demanded quality. However, the initial datum for this solution is already a measure and one may ask, as we did in the introduction, whether there is a measure-valued solution emanating from deterministic initial data which cannot be approximated by a sequence of weak solutions. 

\section{Solution emanating from `atomic' initial data}

We present a non-constant variation of Theorem \ref{thm:constants} proven in \cite{GaWi}:

\begin{theorem}\label{thm.2.1} Let $\Omega\subset\mathbb R^3$ be an open bounded domain, $\mathcal A$ a linear homogeneous constant rank differential operator of order one satisfying $l\geq 3$, and $1\leq p <\infty$. Further, let $\overline z_1, \overline z_2\subset L^\infty(\Omega, \mathbb R^m)$ be such that $\overline z_2 - \overline z_1\notin \Lambda$ a.e.\ in $\Omega$. Assume $z_n:\Omega \mapsto \mathbb R^m$ is an equi-integrable family of functions such that 
\begin{equation*}
\begin{split}
\|z_n\|_p &\leq c <\infty\\
\mathcal A z_n &\to 0 \ \mbox{in}\ W^{-1,r}(\Omega),
\end{split}
\end{equation*}
for some $r\in \left(1,\frac{N}{N-1}\right)$, and $\{z_n\}$ generates a compactly supported Young measure $\nu\in L^\infty_w (\Omega, \mathcal M^1(\mathbb R^m))$ such that
$$
\supp (\nu_x) \subset \left\{  \lambda \overline z_1(x) + (1-\lambda)\overline z_2(x),\ \lambda \in [0,1]\right\}
$$
for a.e.\ $x\in \Omega$.

Then, for a.e. $x\in \Omega$ it holds that
$$
\nu_x  = \delta_{w(x)}
$$
with $w\in L^1(\Omega)$ and $z_n\to w$ in $L^1(\Omega)$.
\end{theorem}

The question of the existence of the demanded measure-valued solution reduces to the question whether there are two weak solutions to \eqref{f.e.cons} $\overline z_1$ and $\overline z_2$ emanating from the same ``atomic" initial conditions such that $\overline z_1-\overline z_2\notin \Lambda$ on a subset of $\Omega$ of positive measure. As noted before, we need to compute the rank of a certain matrix. 

Let $\overline z_1 = (\rho^\alpha,m^\alpha_1, m^\alpha_2, U^\alpha_{11}, U^\alpha_{12}, E^\alpha, r^\alpha_1, r^\alpha_2)$ and $\overline z_2 = (\rho^\beta,m^\beta_1, m^\beta_2, U^\beta_{11}, U^\beta_{12}, E^\beta, r^\beta_1, r^\beta_2)$ be two solutions for which is the constraint \eqref{constrain} effective almost everywhere. The appropriate $Z_{\mathcal A}$ is of the form
\begin{equation}\label{matrix.Z}
\left(
\begin{matrix}
\rho^\alpha - \rho^\beta & m_1^\alpha - m_1^\beta & m_2^\alpha - m_2^\beta\\
m_1^\alpha - m_1^\beta & U_{11}^\alpha - U_{11}^\beta + E^\alpha- E^\beta & U_{12}^\alpha - U_{12}^\beta\\
m_2^\alpha - m_2^\beta & U_{12}^\alpha - U_{12}^\beta & U_{11}^\beta - U_{11}^\alpha  + E^\alpha - E^\beta\\
E^\alpha - E^\beta & r_1^\alpha - r_1^\beta & r_2^\alpha - r_2^\beta
\end{matrix}
\right).
\end{equation}

Below we show the existence of two solutions $(\rho^\alpha, \vv^\alpha, p^\alpha)$ and $(\rho^\beta,\vv^\beta, p^\beta)$ for which the appropriate matrix $Z_{\mathcal A}$ is of rank $3$ on a set of positive measure.

\subsection{The self-similar solution}
\begin{figure}[h!]
  \begin{center}
    \begin{tikzpicture}
		  \draw [thick] (-8,3) to (8,3) node [black,right] {$x_1$};
			\draw (0,3) to (-3,7.5) node [above] {$x_1=st$};
			\node at (-4,5) {$(\rho_-,(v_K,0), p_-)$};
			\node at (3,5) {$(\rho_K, (0,0),p_+)$};
    \end{tikzpicture}
		\caption{Self-similar solution}
		\label{fig.1}
  \end{center}
\end{figure}
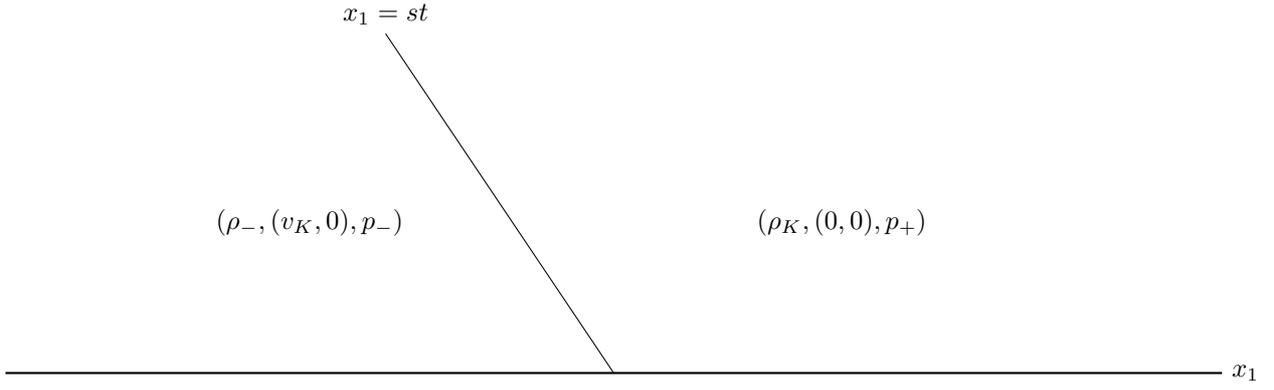

Take Riemann initial data of the following form:
\begin{equation}
(\rho,\vv,p)_0 = \left\{
\begin{array}{l}
(\rho_-,(v_K, 0), p_-)\ \mbox{for }x_1<0\\
(\rho_K,(0,0),p_+)\ \mbox{for }x_1>0,
\end{array}
\right.
\end{equation}
where
$$
\rho_K = \rho_-\frac{p_- + 3p_+}{3p_- + p_+},\quad v_K = \frac{\sqrt 2}{\sqrt \rho_-}\frac{p_+ - p_-}{\sqrt{p_- + 3p_+}}
$$
and $\rho_-,\ p_-,\ p_+ >0$, $p_+>p_-$ are given constants. According to \cite{smoller} there is a self-similar solution consisting of a 1-shock. In particular, let 
$$
s = -\frac{p_+ + 3p_-}{\sqrt{2\rho_-(p_- + 3p_+)}}.
$$
Then a triple
$$
(\rho^\alpha,\vv^\alpha,p^\alpha) = \left\{
\begin{array}{l}
(\rho_-,(v_K,0),p_-)\ \mbox{for }x_1<st,\\
(\rho_K,(0,0),p_+)\ \mbox{for }x_1>st
\end{array}
\right.
$$
is a weak solution to \eqref{f.euler}, see Figure \ref{fig.1}.

\subsection{The wild solution}
Here we present the necessary details of the construction from \cite{KlKrMaMa}. First, according to \cite[Section 4.3]{KlKrMaMa}, we define a pressure $p_\delta = p_+ + \delta_p$ and a velocity $v_\delta = \delta_v$ where $\delta_v = \delta_p \sqrt{\frac{2}{\rho_K(4p_+ + 3\delta_p)}}$. Also, set $\rho_\delta=\rho_K\frac{3p_\delta+p_+}{3p_++p_\delta}$.

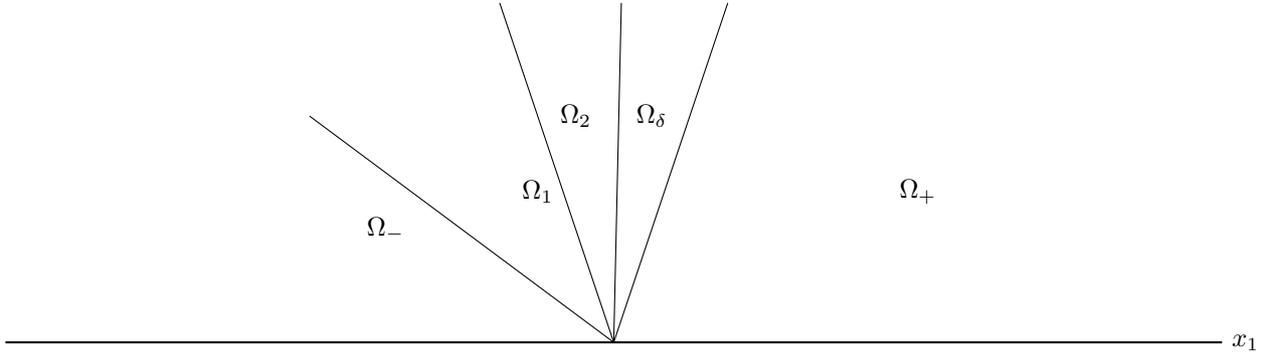
\begin{figure}[h!]
  \begin{center}
	 \begin{tikzpicture}
	   \draw [thick] (-8,3) to (8,3) node [black,right] {$x_1$};
		 \node at (4,5) {$\Omega_+$};
		 \draw  (0,3) to (1.5,7.5);
		 \node at (0.5,6) {$\Omega_\delta$};
		 \draw  (0,3) to (0.1,7.5);
		 \node at (-0.5,6) {$\Omega_2$};
		 \draw  (0,3) to (-1.5,7.5);
		 \node at (-1,5) {$\Omega_1$};
		 \draw  (0,3) to (-4, 6);
		 \node at (-3,4.5) {$\Omega_-$};
	 \end{tikzpicture}
	\caption{Fan partition}
	\label{fig.2}
	\end{center}
\end{figure}

Note that $\delta_v = \delta_v(\delta_p)$ is a smooth function on a neighborhood of $0$, $\delta_v(0) = 0$, and $\delta_p$ is a positive arbitrarily small number.  The time-space is then divided into regions $\Omega_-$, $\Omega_1$, $\Omega_2$, $\Omega_\delta$ and $\Omega_+$ as shown in Figure \ref{fig.2}.

 Between $\Omega_\delta$ and $\Omega_+$ there is a 3-shock. In order to handle the region $\Omega_-\cup \Omega_1\cup \Omega_2\cup\Omega_\delta$, we use the Galilean transformation 
\begin{equation*}
(\rho,(v,u),p)(t,x)\mapsto (\rho,(v-\delta_v,u),p)(t,x+\delta_v t{\bf e}_1)
\end{equation*}
 to get the right state 
$$
(\rho_-, (v_K-\delta_v, 0), p_-) \ \mbox{in }\Omega_-
$$
and the left state
$$(\rho_{\delta_p}, (0,0), p_+ + \delta_p)\ \mbox{in }\Omega_\delta.$$
The Galilean transformation also changes the sets $\Omega_-$, $\Omega_1$, $\Omega_2$ and $\Omega_\delta$. However, we keep the same notation for the sake of simplicity. 

According to \cite[Theorem 3.1]{KlKrMaMa}, there exist infinitely many admissible weak solutions to \eqref{f.euler}. We denote one of these solutions by $(\rho^\beta, (\vv')^\beta,p^\beta)$ -- this solution has to be transformed back by Galilean transformation into the solution $(\rho^\beta, \vv^\beta,p^\beta)$. Such a solution fulfills 
$$
(\rho^\beta,(\vv')^\beta,p^\beta)\restriction_{\Omega_-} = (\rho_-, (v_K-\delta_v,0), p_-),\qquad (\rho^\beta,(\vv')^\beta,p^\beta)\restriction_{\Omega_\delta} = (\rho_\delta, (0,0), p_+ +{\delta_p}).
$$

Moreover, we have the following:
\begin{equation*}
\begin{split}
|(\vv')^\beta\restriction_{\Omega_2}|^2 &= \varepsilon_2 + \tilde\varepsilon_2\\
\rho^\beta\restriction_{\Omega_1} &= \rho_1\in \mathbb R^+\\
\rho^\beta\restriction_{\Omega_2} &= \rho_2\in \mathbb R^+\\
p^\beta\restriction_{\Omega_1} &= p_1\\
p^\beta\restriction_{\Omega_2} &= p_2\\
\end{split}
\end{equation*}
where 
\begin{equation*}
\begin{split}
\varepsilon_2 = \varepsilon_2(v_--\delta_p,\rho_1,\rho_2,p_-,p_+)\\
\tilde\varepsilon_2 = \tilde\varepsilon_2(v_--\delta_p,\rho_1,\rho_2,p_-,p_+)\\
p_2: = p_2(\rho_2,p_-)
\end{split}
\end{equation*}
are continuous functions. Moreover, the construction is such that $\rho_K - \rho_1$ and $\rho_2 - \rho_K$ may be arbitrarily small. As $\delta_v$ and $\delta_p$ are also arbitrarily small, we verify the property that the matrix~\eqref{matrix.Z} has full rank for $\rho_1 = \rho_2 = \rho_K$, $v_- = v_K$ and $p_\delta = p_+$ on a set of positive measure. The aforementioned continuity then allows to use the intended approximation. 

The interface between $\Omega_2$ and $\Omega_\delta$ is $\{x_1 = \delta_vt\}$. Consequently, the domain $\Omega_2\cap \{x_1>st\}$ is nonempty and has positive measure since $s<0$. On this set we take 
\begin{equation*}
\begin{split}
z^\alpha &= \left(\rho^\alpha, \rho^\alpha v^\alpha, \rho^\alpha u^\alpha, \frac 12 \rho^\alpha((v^\alpha)^2 - (u^\alpha)^2), \rho^\alpha u^\alpha v^\alpha, \frac12 \rho^\alpha |\vv^\alpha|^2 + p^\alpha, \left(\frac 12 \rho^\alpha |\vv^\alpha|^2 + 2p^\alpha\right) v^\alpha, \left(\frac 12 \rho^\alpha |\vv^\alpha|^2 + 2p^\alpha \right)u^\alpha\right)\\
z^\beta &= \left(\rho^\beta, \rho^\beta v^\beta, \rho^\beta u^\beta, \frac 12 \rho^\beta((v^\beta)^2 - (u^\beta)^2), \rho^\beta u^\beta v^\beta, \frac12 \rho^\beta |\vv^\beta|^2 + p^\beta, \left(\frac 12 \rho^\beta |\vv^\beta|^2 + 2p^\beta\right) v^\beta, \left(\frac 12 \rho^\beta |\vv^\beta|^2 + 2p^\beta \right)u^\beta\right).
\end{split}
\end{equation*}
Note that $|\vv^\beta|^2 = |(\vv')^\beta + \delta_p{\bf e}_1|^2 = \varepsilon_2 + \tilde \varepsilon_2 + o(\delta)$ on $\Omega_2$. The matrix \eqref{matrix.Z} on the considered set is of the form
\begin{equation*}
\left(
\begin{matrix}
\rho_2 - \rho_K & \rho_2 v^\beta & \rho_2 u^\beta\\
\rho_2 v^\beta & \rho_2 |v^\beta|^2 + p_2 - p_+ & \rho_2 u^\beta v^\beta\\
\rho_2 u^\beta & \rho_2 u^\beta v^\beta & \rho_2 |u^\beta|^2 + p_2 - p_+\\
\frac 12 \rho_2|\vv^\beta|^2 + p_2-p_+ & \left(\frac 12 \rho_2 |\vv^\beta|^2 + 2p_2\right) v^\beta & \left(\frac 12 \rho_2 |\vv^\beta|^2 + 2p_2\right) u^\beta
\end{matrix}
\right).
\end{equation*}

The determinant of the submatrix consisting of the first, second, and third rows is 
$$
(p_2 - p_+) \left((\rho_2 - \rho_K)(p_2 - p_+) - \rho_2\rho_K |\vv^\beta|^2\right).
$$
So the corresponding matrix $Z_{\mathcal A}$ is of rank 3 once we know that $p_2\neq p_+$ and 
\begin{equation}\label{condition.2}
|\vv^\beta|^2\neq \frac{(\rho_2 - \rho_K)(p_2-p_+)}{\rho_2\rho_K}.
\end{equation}
We have $p_2\approx p_-\left(\frac{\rho_K}{\rho_-}\right)^2\neq p_+$ once we know that $p_+>p_-$, and so it remains to verify \eqref{condition.2}.

Since all above mentioned quantities are bounded and $\rho_2 - \rho_K$ is negligible, the right hand side of \eqref{condition.2} can be made arbitrarily close to zero. We need to show that $|\vv^\beta|^2$ is far away of zero. According to \cite[Sections 3.2 \& 3.8]{KlKrMaMa} we have
$$
|\vv^\beta|^2 \approx \varepsilon_2(v_K,\rho_K,\rho_K,p_-,p_+) + \tilde\varepsilon_2(v_K,\rho_K,\rho_K,p_-,p_+)= 4\frac{(p_+ - p_-)(p_++p_-)^2}{\rho_-(3p_+ + p_-)(3p_-+p_+)}.
$$
Consequently, \eqref{condition.2} is fulfilled and it is allowed to consider also the intended small perturbations. The Young measure
$$
\nu_{t,x}  = \frac12\delta_{(\rho^\alpha,\rho^\alpha \vv^\alpha, \frac12\rho^\alpha|\vv^\alpha|^2+p^\alpha)} + \frac12\delta_{(\rho^\beta,\rho^\beta \vv^\beta, \frac12\rho^\beta|\vv^\beta|^2+p^\beta)}
$$
is a measure value solution which, due to Theorem \ref{thm.2.1}, cannot be generated by weak solutions.

We get the following claim as a result of the previous considerations.
\begin{theorem}\label{mainthm}
There exists a measure-valued solution to \eqref{f.euler} with non-constant entropy, emanating from certain Riemann initial data, which cannot be generated by a sequence of weak solutions.
\end{theorem}

{\bf Acknowledgement:} The research of V.M. was supported by the Czech Science Foundation, Grant Agreement GA18--05974S, in the framework of RVO:67985840.

\end{document}